\documentclass[11pt]{article}
\def\bn{\hbox{\it I\hskip -2pt N}}
\def\bz{\hbox{\it Z\hskip -4pt Z}}
\def\bq{\hbox{\it l\hskip -5.5pt Q}}

\def\demo{\noindent{\bf Proof.}}
\newtheorem{theorem}{Theorem}[section]
\newtheorem{lemma}[theorem]{Lemma}
\newtheorem{remark}[theorem]{Remark}
\newtheorem{example}[theorem]{Example}
\newtheorem{Corollary}[theorem]{Corollary}
\newtheorem{definition}[theorem]{Definition}
\newtheorem{proposition}[theorem]{Proposition}
\begin{document}
\begin{center}
{\uppercase{\bf Projections of cones and the arithmetical rank of toric varieties}}
\end{center}
\advance\baselineskip-3pt
\vspace{2\baselineskip}
\begin{center}
{{\sc Anargyros Katsabekis }\\
\vspace{2\baselineskip}
\vspace{\baselineskip}}
\end{center}
\vskip.5truecm\noindent
\vskip.5truecm\noindent\vspace*{5pt}
\begin{center}
{\uppercase{abstract}}
\end{center}

Let $I_M$ and $I_N$ be defining ideals of toric varieties such that $I_M$ is a projection of $I_N$, i.e. $I_N \subseteq I_M$. We give necessary and sufficient conditions for the equality $I_M=rad(I_N+(f_1,\ldots,f_s))$, where $f_1,\ldots,f_s$ belong to $I_M$. Also a method for finding toric varieties which are set-theoretic complete intersection is given. Finally we apply our method in the computation of the arithmetical rank of certain toric varieties and provide the defining equations of the above toric varieties.

\section{Introduction} One of the classical problems of Algebraic Geometry with a long history, see \cite{Eis} Chapter 15, is to determine the minimal number of equations needed to define set-theoretically an algebraic variety over an algebraically closed field. Even more difficult is to provide minimal sets of equations that define the algebraic variety. The problem is open even for very simple cases, like the Macaulay curve $(t^4,t^3u,tu^3,u^4)$ in the three dimensional projective space. This article addresses these two problems for toric varieties and in several cases we are able to compute the minimal number, but also provide the equations that define the variety set-theoretically.\\
\par Let $M=(a_{i,j})$ be an $m \times n$ matrix with integer entries $a_{i,j}$, such that every column has atleast one non-zero entry. Let $K$ be a field and let $A=\{{\bf a}_1,\ldots,{\bf a}_n \}$ be the set of vectors in $\bz ^m$, where ${\bf a}_i=(a_{1,i},\ldots,a_{m,i})$ for $1 \leq i \leq n$. The {\em toric ideal} $I_M$ associated with $M$ is the kernel of the $K$-algebra homomorphism
$$\phi: K[x_1,\ldots, x_n]\rightarrow
 K[t_1,\ldots,t_{m},t_1^{-1},\ldots ,t_{m}^{-1}]$$ given by $$\phi(x_i)={\bf t}^{{\bf a}_i} :=t_1^{a_{1,i}} \cdots t_m^{a_{m,i}} \qquad \mbox{for  
all }i=1,\ldots,n.$$ The ideal $I_M$ is prime and therefore $rad(I_M)=I_M$. A difference of two monomials is called a {\em binomial}. Every vector ${\bf u}$ in $\bz^n$ can be written uniquely as ${\bf u}={\bf u}_+-{\bf u}_-$, where ${\bf u}_+=(u_{+,1},\ldots,u_{+,n})$ and ${\bf u}_- =(u_{-,1},\ldots,u_{-,n})$ are non-negative and have disjoint support. If $D$ is an $m \times n$ matrix with rational entries and columns $\{ {\bf d}_1,\ldots,{\bf d}_n \}$, then the kernel of $D$ is $$ker(D)=\{ {\bf u}=(u_1,\ldots,u_n) \in \bq^n | u_1{\bf d}_1+\cdots+u_n{\bf d}_n={\bf 0} \}.$$ Set $ker_{\bz}(D)=ker(D) \cap \bz^n$. The height $ht(I_M)$ of $I_M$ equals the rank of the lattice $ker_{\bz}(M)$ (see \cite{St}).
\begin{lemma} (\cite{St}) The toric ideal $I_M$ is generated by the binomials ${\bf x}^{{\bf u}_+}-{\bf x}^{{\bf u}_-}$, where ${\bf u}$ belongs to $ker_{\bz}(M)$.\\ 
\end{lemma}
\par We grade the polynomial ring $K[x_1,\ldots,x_n]$ by setting $deg_A(x_i)={\bf a}_i$ for $i=1,\ldots,n$. We define the $A$-{\em degree} of the monomial ${\bf x}^{\bf u}$ to be $$deg_A({\bf x}^{\bf u}) :=u_1{\bf a}_1+\cdots+u_n{\bf a}_n \in \bn(A),$$ where $\bn(A)$ is the semigroup generated by $A$. Every toric ideal $I_M$ is $A$-homogeneous, since it is generated by binomials and every binomial ${\bf x}^{{\bf u}_+}-{\bf x}^{{\bf u}_-}$ is $A$-homogeneous.\\
\par The {\em toric variety} $X_M$ associated with $M$ is the set $V(I_M) \subset K^n$ of zeroes of $I_M$ in the sense of \cite{St}, which also includes non-normal varieties. The {\em toric set} $\Gamma(M)$ determined by $M$ is the subset of $K^n$ defined parametrically by $x_i=t_1^{a_{1,i}} \cdots t_m^{a_{m,i}}$ for all $i$, i.e. it is the set of points that can be expressed in the form $$(t_1^{a_{1,1}}\cdots t_m^{a_{m,1}}, \ldots,t_1^{a_{1,i}} \cdots t_m^{a_{m,i}},\ldots,t_1^{a_{1,n}}\cdots t_m^{a_{m,n}})$$ for some $t_i$ in $K$. Note that $\Gamma(M)$ is a subset of $X_M$. When $m=1$ and $a_{1,1} < a_{1,2} < \cdots < a_{1,n}$ are positive integers, the g.c.d. of which equals 1, then $\Gamma(M)$ is known as a monomial curve and $I_M$ as the ideal of the monomial curve.\\
\par We associate to the toric variety $X_M$ the rational polyhedral cone $\sigma=pos_{\bq}(A):=\{\sum_{i=1}^{n}d_{i}{\bf a}_i|d_i \in \bq \ and \ d_i \geq 0 \}$. The dimension of $\sigma$ is equal to the dimension of the vector space ${\bq}A=\{\sum_{i=1}^{n}d_{i}{\bf a}_i|d_i \in \bq \}$ and also is equal to the dimension of $X_M$.\\
\par In this paper we consider the following two problems related to the toric ideal $I_M$: 

\noindent (I) Given a toric ideal $I_N$, such that $I_N \subseteq I_M$, and a set of binomials $\{ f_1,\ldots,f_s \}$ in $I_M$, formulate a criterion for deciding the equality $I_M=rad(I_N+(f_1,\ldots,f_s))$. \\
\noindent (II) Find the smallest number of polynomials needed to generate $I_M$ up to radical. This problem is more general from the corresponding problem in Algebraic Geometry of the determination of the minimum number of equations needed to define a toric variety $X_M$ set-theoretically, over an algebraically closed field. This number is called {\em arithmetical rank} of $X_M$ and will be denoted by $ara(X_M)$. The Generalized Krull's principal ideal theorem provides a lower bound for the arithmetical rank of $X_M$, namely the height of $I_M$. When $ht(I_M)=ara(X_M)$ the ideal $I_M$ (and the variety $X_M$ as well) is called a {\em set-theoretic complete intersection}.\\

\par The problem (I) was studied by Eliahou-Villarreal in \cite{E-V} in the special case that $I_N=(0)$. There they give a nessessary and sufficient condition for the equality $rad(f_1,\ldots,f_s)=I_M$. More precisely they prove that:
\begin{theorem} (\cite{E-V},\cite{Vil}) Let $\{ f_1,\ldots,f_s \}$ be a set of binomials in the toric ideal $I_M$. Set $J=(f_1,\ldots,f_s)$ and $G=<{\widehat{f_1}},\ldots,{\widehat{f_s}}> \subset ker_{\bz}(M)$, where for a binomial $f={\bf x}^{\bf u}-{\bf x}^{\bf v} \in K[x_1,\ldots,x_n]$ we let $\widehat{f}={\bf u}-{\bf v} \in \bz^n$. If $char(k)=p \neq 0$ (resp. $char(k)=0$), then $rad(J)=I_M$ if and only if:\\
\noindent (a) $p^k ker_{\bz}(M) \subset G$ for some $k \in \bn$ (resp. $ker_{\bz}(M)=G$),\\
\noindent (b) $rad(J,x_i)=rad(I_M,x_i)$, for all $i=1,\ldots,n$.
\end{theorem}
We generalize (see Theorem 4.4) this result in terms of projections of ideals. Our criterion can be used also to determine different binomial generators for the radical of the ideal of a toric variety (see example 4.6). However in the case that we can make a good choice of a projection $I_N$, minimal binomial generators up to radical for the ideal $I_M$ of the toric variety are derived.\\
\par Basic ingredient of our approach is the notion of projections of toric ideals. This notion, although it was never before explicitly defined, has been used for the first time by J. Herzog in \cite{Her} to prove that the ideal of the monomial curve $(t^{a_1},t^{a_2},t^{a_3})$ is set-theoretic complete intersection. In \cite{Th} the same notion has been used to prove that, when $char(K)=0$, smooth monomial curves are not binomial set-theoretic complete intersections, except for the twisted cubic. Finally in \cite{Tho} A. Thoma used this notion to deduce that certain ideals of monomial curves $(t^{a_1},\ldots,t^{a_n})$ are set-theoretic complete intersections. The techniques developed there can not always be applied, for example the ideal of the monomial curve $(t^4,t^6,t^{11},t^{13})$ was generally unknown whether or not it is a set-theoretic complete intersection. Our method begins with a toric ideal which is set-theoretic complete intersection and it produces a large number of toric ideals which are set-theoretic complete intersections. This method also provides the defining equations of the toric variety. \\ 
In section 2 we introduce the basic notion of this paper, the notion of projection of toric ideals, and present its connection with the geometric notion of projection of cones.\\
In section 3 we give  necessary and sufficient conditions for the equality $rad(I_N+(f_1,\ldots,f_s))=I_M$, where $I_M$ is a projection of $I_N$ and $f_1,\ldots,f_s$ belong to $I_M$.\\
In section 4 we study the previous equality in the special case that $f_1,\ldots,f_s$ are binomials.\\
In section 5 we develop a method for finding toric ideals which are set-theoretic complete intersections.\\
In section 6 we apply the theory developed in section 5 in the computation of the exact value of the arithmetical rank of certain toric ideals. Among other results, we prove that the ideal of the monomial curve $(t^4,t^6,t^{a},t^{b})$ is set-theoretic complete intersection, so for $a=11, b=13$ the ideal of the monomial curve $(t^4,t^6,t^{11},t^{13})$ is set-theoretic complete intersection.
\section{Projections of toric ideals} We consider the toric ideals $I_M$, $I_N$ associated with the  $m \times n$ matrix $M=(a_{i,j})$ and $l \times n$ matrix $N=(b_{i,j})$ respectively. Let $A=\{{\bf a}_1,\ldots,{\bf a}_n \}$, $B=\{{\bf b}_1,\ldots,{\bf b}_n \}$, where ${\bf a}_i=(a_{1,i},\ldots,a_{m,i})$ and ${\bf b}_i=(b_{1,i},\ldots,b_{l,i})$ for $1 \leq i \leq n$.
\begin{definition} We say that $I_M$ is a {\em projection} of $I_N$ if $I_N \subseteq I_M$.
\end{definition}
Let $\widehat{\pi}: {\bq}^l \rightarrow {\bq}^m$ be a rational affine map with $\widehat{\pi}(pos_{\bq}(B))=pos_{\bq}(A)$. We call $$\pi:=\widehat{\pi}|_{pos_{\bq}(B)}: pos_{\bq}(B) \rightarrow pos_{\bq}(A)$$ a {\em projection of cones}.\\ 
\par The next Theorem makes the connection between the algebraic notion of projection of toric ideals and the geometric notion of projections of cones.  
\begin{theorem} The following are equivalent:\\
\noindent (a) $I_M$ is a projection of $I_N$.\\
\noindent (b) The lattice $ker_{\bz}(N)$ is a subset of the lattice $ker_{\bz}(M)$.\\
\noindent (c) Every $B$-homogeneous ideal in $K[x_1,\ldots,x_n]$ is also $A$-homogeneous.\\
\noindent (d) There is a projection of cones $\pi: pos_{\bq}(B) \rightarrow pos_{\bq}(A)$ given by $\pi({\bf b}_i)={\bf a}_i$ for all $i=1,\ldots,n$.\\
\noindent (e) There is an $m \times l$ matrix $D$ with rational entries such that $DN=M$.
\end{theorem}
\demo (a) $\Rightarrow$ (b) Let ${\bf u}={\bf u}_{+}-{\bf u}_{-}$ be an element of $ker_{\bz}(N)$. Then the binomial ${\bf x}^{{\bf u}_+}-{\bf x}^{{\bf u}_-}$ belongs to $I_N$ which, from the assumption, is a subset of $I_M$. Consequently ${\bf u}$ is in $ker_{\bz}(M)$.\\
\noindent (b) $\Rightarrow$ (c) Let $I \subseteq K[x_1,\ldots x_n]$ be a $B$-homogeneous ideal and ${\bf x}^{\bf u}$, ${\bf x}^{\bf v}$ two monomials of a $B$-homogeneous generator $f$ of $I$. Where ${\bf u}=(u_1,\ldots,u_{n})$ and ${\bf v}=(v_1,\ldots,v_n)$. Set $g={\bf x}^{\bf u}-{\bf x}^{\bf v}$. We have $$u_1{\bf b}_1+ \cdots +u_n{\bf b}_n=v_1{\bf b}_1+ \cdots +v_n{\bf b}_n,$$ which implies that ${\widehat{g}}$ belongs to the lattice $ker_{\bz}(N)$. As $ker_{\bz}(N) \subseteq ker_{\bz}(M)$, we obtain that the vector ${\widehat{g}}$ belongs to $ker_{\bz}(M)$ and so $$u_1{\bf a}_1+ \cdots +u_n{\bf a}_n=v_1{\bf a}_1+ \cdots +v_n{\bf a}_n.$$ Hence $I$ is $A$-homogeneous.\\
\noindent (c) $\Rightarrow$ (d) It is enough to consider the case $dim(pos_{\bq}(B))=l$. Let $\{ {\bf b}_{j_1},\ldots,{\bf b}_{j_l} \}$ be a base for the ${\bq}$-vector space ${\bq}^l$ and define $\widehat{\pi}: {\bq}^l \rightarrow {\bq}^m$ by $\widehat{\pi}({\bf b}_{j_i})={\bf a}_{j_i}$ for all $i=1,\ldots,l$. Set $\pi=\widehat{\pi}|_{pos_{\bq}(B)}$. Obviously $\pi({\bf b}_{j_i})={\bf a}_{j_i}$ for every $i=1,\ldots,l$. Let ${\bf b}_{j_i}=\sum_{k=1}^{l} \lambda_{k}{\bf b}_{j_k} \in pos_{\bq}(B)$ for some $\lambda_{k} \in \bq$ and $i \in \{ l+1,\ldots,n \}$. Clear the denominators to obtain an equality ${\nu}{\bf b}_{j_i}=\sum_{k=1}^{l} \xi_{k} {\bf b}_{j_k}$, where ${\nu}$ and ${\xi}_1,\ldots,{\xi}_l$ are integers. Suppose that ${\nu}$ is positive, the case ${\nu}<0$ is essentially the same. The ideal $I=(x_{j_1}^{{\xi}_{1,-}} \cdots x_{j_l}^{{\xi}_{l,-}}x_{j_i}^{\nu}-x_{j_1}^{{\xi}_{1,+}} \cdots x_{j_l}^{{\xi}_{l,+}})$ is $B$-homogeneous and therefore $A$-homogeneous. Thus ${\nu}{\bf a}_{j_i}=\sum_{k=1}^{l} {\xi}_{k} {\bf a}_{j_k}$, which leads to ${\bf a}_{j_i}=\sum_{k=1}^{l} \lambda_{k}{\bf a}_{j_k}$. Hence $\pi({\bf b}_{j_i})={\bf a}_{j_i}$, for every $i \in \{ 1,\ldots,n \}$.\\
\noindent (d) $\Rightarrow$ (e) The matrix $D$ is the matrix of $\widehat{\pi}$ in the canonical bases of ${\bq}^l$ and ${\bq}^m$.\\
\noindent (e) $\Rightarrow$ (b) Let ${\bf u}=(u_1,\ldots,u_n)$ be an element of $ker_{\bz}(N)$. Then $D(N{\bf u}^{\bf T})={\bf 0}^{\bf T}$, and therefore $M{\bf u}^{\bf T}={\bf 0}^{\bf T}$. So $ker_{\bz}(N) \subseteq ker_{\bz}(M)$.\\
\noindent (b) $\Rightarrow$ (a) Let $f={\bf x}^{{\bf u}_+}-{\bf x}^{{\bf u}_-}$ be a binomial generator of $I_N$, where ${\widehat{f}}$ belongs to $ker_{\bz}(N)$. As $ker_{\bz}(N) \subseteq ker_{\bz}(M)$, we take that ${\widehat{f}}$ is in $ker_{\bz}(M)$. Hence $f$ belongs to $I_M$. $\diamondsuit $
\begin{Corollary} If $I_M$ is a projection of $I_N$, then $ht(I_N) \leq ht(I_M)$.
\end{Corollary}
\demo We have $ker_{\bz}(N) \subseteq ker_{\bz}(M)$ and therefore $rank(ker_{\bz}(N)) \leq rank(ker_{\bz}(M))$, by Theorem 2.2. Thus $ht(I_N)=rank(ker_{\bz}(N)) \leq rank(ker_{\bz}(M))=ht(I_M)$. $\diamondsuit $
\section{Set-Theoretic Generation} 
Let $I_M$, $I_N$ be toric ideals such that $I_M$ is a projection of $I_N$. Where $M$ is an $m \times n$ integer matrix and $N=(b_{i,j})$ is an $l \times n$ matrix with non-negative integer entries. Set $B=\{ {\bf b}_1,\ldots,{\bf b}_n \}$, where ${\bf b}_i=(b_{1,i},\ldots,b_{l,i})$ for $i=1,\ldots,n$. Note that if $pos_{\bq}(B)$ is strongly convex, i.e. ${\bf 0}$ is the only invertible element of ${\bn}(B)$, then we can choose an appropriate matrix $N$ with non-negative integer entries. The toric ideal $I_N$ is the kernel of the $K$-algebra homomorphism $$\phi: K[x_1,\ldots,x_n] \rightarrow K[t_1,\ldots,t_l]$$ given by $$\phi(x_i)={\bf t}^{{\bf b}_i} \qquad \mbox{for all } i=1,\ldots,n.$$
We shall denote by $(\phi(I_M))^{e}$ the ideal $\phi(I_M)K[t_1,\ldots,t_l]$.
\begin{theorem} Assume that $\Gamma(N)=V(I_N)$ in ${\bar{K}}^n$. Let $\{ f_1,\ldots,f_s \}$ be a set of polynomials in $I_M$ and let $J=(f_1,\ldots,f_s)$. Then $$I_{M}=rad(I_{N}+J)$$ if and only if $$rad((\phi(I_{M}))^{e})=rad(\phi(f_1),\ldots,\phi(f_s)).$$
\end{theorem}
\demo Let us first assume that   $$rad((\phi(I_{M}))^{e})=rad(\phi(f_1),\ldots,\phi(f_s)).$$ We claim that $$I_{M}=rad(I_{N}+J).$$ By Hilbert's Nullstellensatz it is enough to prove that any point ${\bf x}=(x_1,\ldots,x_n)$ of $V(I_N,J)$ in ${\bar{K}}^n$ belongs to $V(I_{M})$. This ${\bf x}$ belongs also to $V(I_{N})=\Gamma(N)$. The last statement means that there exist $T_i \in \bar{K}$ for $1 \leq i \leq l$ such that $x_i={\bf T}^{{\bf b}_i}$. Note that $f({\bf x})=\phi(f)(T_1,\ldots,T_l)$ for every $f \in K[x_1,\ldots,x_n]$, since ${\bf x}$ is of the above form. In addition the point $(T_1,\ldots,T_l)$ belongs to $V(\phi(f_1),\ldots,\phi(f_s))=V((\phi(I_M))^{e})$, because $f_i({\bf x})=0$ for every $i \in \{1,\ldots,s \}$. Let $f$ be a polynomial in $I_{M}$. Then ${\phi}(f)$ belongs to $(\phi(I_{M}))^{e}$ and therefore ${\phi}(f)(T_1,\ldots,T_l)=0$. Thus $f({\bf x})=0$, which implies the demanded relation. 
Conversely assume that $I_{M}=rad(I_{N}+J)$. We will prove that $$rad((\phi(I_{M}))^{e})=rad((\phi(I_{N}))^{e}+(\phi(f_1),\ldots,\phi(f_s))).$$ It is enough to prove it for a generator $f=\phi(g)$ of $(\phi(I_{M}))^{e}$, where $g \in I_{M}$. We have $g^k=h_1+h_2$ for some $k \in \bn$, where $h_1 \in I_{N}$ and $h_2 \in J$. Therefore $f^k=\phi(g^k)=\phi(h_1)+\phi(h_2)$, which means that $f$ belongs to $rad((\phi(I_{N}))^{e}+(\phi(f_1),\ldots,\phi(f_s)))$. But $I_{N}=ker(\phi)$, so $rad((\phi(I_{M}))^{e})=rad(\phi(f_1),\ldots,\phi(f_s))$. $\diamondsuit $
\begin{remark} {\rm Every toric variety can always be expressed as an appropriate toric set over an algebraically closed field, for details see \cite{K-T}. The proof of this fact is constructive and also an algorithm is given there to find this toric set. Therefore, in any case, the condition $\Gamma(N)=V(I_{N})$ in ${\bar{K}}^n$ can be achieved, by choosing an appropriate matrix $N$ for any toric ideal.}
\end{remark}
\begin{remark} {\rm We can not omit from the assumptions of Theorem 3.1 the fact that $V(I_{N})$ coincides with the toric set $\Gamma(N)$. Let $I_N$ be the toric ideal associated with the matrix $$N=\left(\matrix{2 & 1 & 0 & 0 & 1 & 2 \cr
1 & 2 & 2 & 1 & 0 & 0 \cr
0 & 0 & 1 & 2 & 2 & 1 \cr} \right),$$ and let $I_M$ the toric ideal associated with the matrix $$M=\left(\matrix{3 & 3 & 2 & 1 & 1 & 2 \cr
2 & 1 & 1 & 2 & 3 & 3 \cr} \right).$$ Notice that $I_M$ is a projection of $I_N$ and that the toric variety $V(I_N)$ does not coincide with $\Gamma(N)$(see \cite{K-T}). The toric ideal $I_M \subset K[x_1,\ldots,x_6]$ is minimally generated by the following $12$ binomials : $x_1^2-x_2x_3x_4,x_3^3-x_1x_2,x_4^3-x_5x_6,x_6^2-x_3x_4x_5,{x_3^2}x_5-x_2{x_4^2},x_3{x_5^2}-{x_4^2}x_6,x_1{x_3^2}-{x_2^2}x_4,{x_3^2}x_4-x_2x_6,{x_3^2}x_5-x_1x_6,x_3{x_4^2}-x_1x_5,x_2x_5-x_3x_6,x_1x_4-x_3x_6$. Let $\phi: K[x_1,\ldots,x_6] \rightarrow K[t_1,t_2,t_3]$ be the $K$-algebra homomorphism with $I_N=ker(\phi)$. The ideal $({\phi}(I_M))^{e}$ is minimally generated by the binomials : $t_1^{4}t_2^{2}-t_1t_2^{5}t_3^{3},t_2^{6}t_3^{3}-t_1^{3}t_2^{3},t_2^{3}t_3^{6}-t_1^{3}t_3^{3},t_1^{4}t_3^{2}-t_1t_2^{3}t_3^{5},t_2^{5}t_3^{4}-t_1^{3}t_2^{2}t_3,t_1t_2^{4}t_3^{4}-t_1^{4}t_2t_3,t_2^{4}t_3^{5}-t_1^{3}t_2t_3^{2}$. Therefore $rad(({\phi}(I_M))^{e})=rad(t_2^6t_3^3-t_1^3t_2^3,t_2^3t_3^6-t_1^3t_3^3)$, since a power of the other generators of $({\phi}(I_M))^{e}$ belongs to the ideal $(t_2^6t_3^3-t_1^3t_2^3,t_2^3t_3^6-t_1^3t_3^3)$. Observe that $t_2^6t_3^3-t_1^3t_2^3=\phi(x_3^3-x_1x_2)$ and $t_2^3t_3^6-t_1^3t_3^3=\phi(x_4^3-x_5x_6)$. But $I_M \neq rad(I_N+(x_3^3-x_1x_2,x_4^3-x_5x_6))$, since $I_N$ does not have any monic binomial in the variable $x_6$ and therefore no power of the binomial $x_6^2-x_3x_4x_5 \in I_M$ belongs to the ideal $I_N+(x_3^3-x_1x_2,x_4^3-x_5x_6)$.
Let $$D=\left(\matrix{2 & 1 & 0 & 0 & 1 & 2 \cr
1 & 2 & 2 & 1 & 0 & 0 \cr
0 & 0 & 1 & 2 & 2 & 1 \cr
2 & 2 & 1 & 0 & 0 & 1 \cr
0 & 1 & 2 & 2 & 1 & 0 \cr
1 & 0 & 0 & 1 & 2 & 2 \cr} \right).$$ 
Then $\Gamma(D)=V(I_N)$, a proof of this fact can be found in \cite{K-T}. Let ${\psi}:K[x_1,\ldots,x_6] \rightarrow K[t_1,\ldots,t_6]$ be the $K$-algebra homomorphism with $Ker({\psi})=I_D=I_N$. We have $rad(({\psi}(I_M))^{e}) \neq rad({\psi}(x_3^3-x_1x_2),{\psi}(x_4^3-x_5x_6))$, since no power of $\psi(x_6^2-x_3x_4x_5)$ belongs to the ideal $(\psi(x_3^3-x_1x_2),\psi(x_4^3-x_5x_6))$.}
\end{remark}
\begin{definition} The toric ideal $I_M$ is called set-theoretic complete intersection on $I_N$ if there are polynomials $f_1,\ldots,f_s$ in $I_M$, where $s$ is equal to the difference of the heights of $I_M$ and $I_N$, satisfying $I_M=rad(I_N+(f_1,\ldots,f_s))$.\\
\end{definition}
\par The next Corollary is directly derived from Theorem 3.1.
\begin{Corollary} Keep the assumptions of Theorem 3.1. The toric ideal $I_{M}$ is set-theoretic complete intersection on $I_{N}$ if and only if the radical of the ideal $(\phi(I_{M}))^{e}$ is equal to the radical of an ideal generated by $ht(I_M)-ht(I_N)$ elements of the form $\phi(f)$, where $f \in I_{M}$.
Moreover, if $I_N$ is set-theoretic complete intersection and $I_M$ is set-theoretic complete intersection on $I_N$ then $I_M$ is set-theoretic complete intersection.
\end{Corollary}
\section{Binomial generation} In this section we will give an equivalent condition for the equality $I_{M}=rad(I_{N}+(f_1,\ldots,f_s))$ when $f_1,\ldots,f_s$ are binomials in $I_M$.\\
\par We consider the toric ideals $I_M$ and $I_N$, where $M$ is an $m \times n$ integer matrix and $N=(b_{i,j})$ is an $l \times n$ integer matrix with non-negative entries. Assume that $I_M$ is a projection of $I_N$ and let $B=\{{\bf b}_1,\ldots,{\bf b}_n \}$, where ${\bf b}_i=(b_{1,i},\ldots,b_{l,i})$ for $i=1,\ldots,n$. The toric ideal $I_N$ is the kernel of the $K$-algebra homomorphism $$\phi:K[x_1,\ldots,x_n] \rightarrow K[t_1,\ldots,t_l]$$ defined by $$\phi(x_i)={\bf t}^{{\bf b}_i} \qquad \mbox{for all } i=1,\ldots,n.$$\\
\par Given a lattice $L \subset \bz^l$ the ideal $$I_{L}:=(\{{\bf t}^{{\bf z}_+}-{\bf t}^{{\bf z}_-}| {\bf z}={\bf z}_{+}-{\bf z}_{-} \in L \}) \subset K[t_1,\ldots,t_l]$$ is called {\em lattice ideal}.
The height of $I_{L}$ equals the rank of the lattice $L$ (see \cite{E-S}). For a prime number $p$ we denote by $L:p^{\infty}$ the lattice $$\{{\bf u} \in \bz^l | p^k{\bf u} \in L \ \mbox{ for some } \ k \in \bn \}.$$ Let ${\bz}B=\{{\bf y}_{\bf u}:=u_1{\bf b}_1+\cdots+u_n{\bf b}_n | {\bf u}=(u_1,\ldots,u_n) \in \bz^n \}$ be the lattice spanned by $B$ and let $$N(ker_{\bz}(M))=\{{\bf y}_{\bf u}=u_1{\bf b}_1+\cdots+u_n{\bf b}_n | {\bf u}=(u_1,\ldots,u_n) \in ker_{\bz}(M) \} \subset {\bz}B.$$ Note that if $ker_{\bz}(M)=<{\bf u}_1,\ldots,{\bf u}_r>$, then $N(ker_{\bz}(M))=<{\bf y}_{{\bf u}_1},\ldots,{\bf y}_{{\bf u}_r}>$.
\begin{lemma} The lattice ideal $I_{N(ker_{\bz}(M))}$ coincides with the ideal $(\phi(I_M))^{e}:(t_1 \cdots t_l)^{\infty}$.
\end{lemma}
\demo Let $f=x_1^{u_1} \cdots x_n^{u_n}-x_1^{v_1} \cdots x_n^{v_n}$ be a binomial in $I_M$, where ${\widehat{f}}=(u_1-v_1,\ldots,u_n-v_n) \in ker_{\bz}(M)$. Then $\phi(f)={\bf t}^{u_1{\bf b}_1} \cdots {\bf t}^{u_n{\bf b}_n}-{\bf t}^{v_1{\bf b}_1} \cdots {\bf t}^{v_n{\bf b}_n}$. Set ${\bf z}=u_1{\bf b}_1+ \cdots +u_n{\bf b}_n-({\bf y}_{\widehat{f}})_+=v_1{\bf b}_1+ \cdots +v_n{\bf b}_n-({\bf y}_{\widehat{f}})_-$ and observe that ${\bf z} \in \bn^l$. We have $\phi(f)={\bf t}^{\bf z}({\bf t}^{({\bf y}_{\widehat{f}})_+}-{\bf t}^{({\bf y}_{\widehat{f}})_-})$ and therefore $(\phi(I_M))^{e} \subseteq I_{N(ker_{\bz}(M))} \subseteq (\phi(I_M))^{e}:(t_1 \cdots t_l)^{\infty}$. Clearly $(\phi(I_M))^{e}:(t_1 \cdots t_l)^{\infty} \subseteq I_{N(ker_{\bz}(M))}:(t_1 \cdots t_l)^{\infty} \subseteq ((\phi(I_M))^{e}:(t_1 \cdots t_l)^{\infty}):(t_1 \cdots t_l)^{\infty}$. But $I_{N(ker_{\bz}(M))}=I_{N(ker_{\bz}(M))}:(t_1 \cdots t_l)^{\infty}$ (see \cite{E-S}) and $(\phi(I_M))^{e}:(t_1 \cdots t_l)^{\infty}=((\phi(I_M))^{e}:(t_1 \cdots t_l)^{\infty}):(t_1 \cdots t_l)^{\infty}$. Thus $I_{N(ker_{\bz}(M))}=(\phi(I_M))^{e}:(t_1 \cdots t_l)^{\infty}$. $\diamondsuit $
\begin{proposition} Let $G$ be a sublattice of $N(ker_{\bz}(M))$. If $char(K)=p \neq 0$ (resp. $char(K)=0$), then the following two conditions are equivalent:\\
\noindent (a) $rad(I_{N(ker_{\bz}(M))})=rad(I_G)$,\\
\noindent (b) $p^k N(ker_{\bz}(M)) \subset G$ for some $k \in \bn$ (resp. $N(ker_{\bz}(M))=G$).
\end{proposition}
\demo Suppose first that $rad(I_{N(ker_{\bz}(M))})=rad(I_G)$. By \cite{E-S} Corollary 2.2 it follows that in characteristic zero $I_{N(ker_{\bz}(M))}=I_G$ and so $N(ker_{\bz}(M))=G$, since a binomial $f={\bf t}^{\bf u}-{\bf t}^{\bf v}$ lies in a lattice ideal $I_L$ if and only if ${\widehat{f}}$ belongs to $L$. Also in characteristic $p \neq 0$ it holds $I_{N(ker_{\bz}(M)):p^{\infty}}=I_{G:p^{\infty}}$ and so $N(ker_{\bz}(M)):p^{\infty}=G:p^{\infty}$. Note that in \cite{E-S} the lattice $L:p^{\infty}$ is denoted by $Sat_{p}(L)$. Suppose that $N(ker_{\bz}(M))=<{\bf y}_{{\bf u}_1},\ldots,{\bf y}_{{\bf u}_r}>$. We have $N(ker_{\bz}(M)) \subset N(ker_{\bz}(M)):p^{\infty}=G:p^{\infty}$. Hence for every $i=1,\ldots,r$ there exist $k_i \in \bn$ such that $p^{k_i}{\bf y}_{{\bf u}_i} \in G$. By choosing $k$ the maximum of all $k_i$ we take $p^k N(ker_{\bz}(M)) \subset G$. Conversely it is clear that in characteristic zero $I_{N(ker_{\bz}(M))}=I_G$ and therefore $rad(I_{N(ker_{\bz}(M))})=rad(I_G)$. We restrict now our attention in the case $char(K)=p \neq 0$. It is enough to show that $N(ker_{\bz}(M)):p^{\infty}=G:p^{\infty}$. Obviously $G:p^{\infty} \subset N(ker_{\bz}(M)):p^{\infty}$. Let ${\bf u} \in N(ker_{\bz}(M)):p^{\infty}$. Then there exist $d \in \bn$ such that $p^d{\bf u}$ is in $N(ker_{\bz}(M))$, so from the hypothesis $p^{d+k} {\bf u}$ belongs to $G$. Thus ${\bf u} \in G:p^{\infty}$ and therefore $rad(I_{N(ker_{\bz}(M))})=rad(I_G)$. $\diamondsuit $
\begin{remark} {\rm From the proof of the above proposition we can see that in characteristic $p \neq 0$ condition (b) is equivalent with the condition $N(ker_{\bz}(M)):p^{\infty}=G:p^{\infty}$.}
\end{remark}
\begin{theorem} Assume that $\Gamma(N)=V(I_N)$ in ${\bar{K}}^n$. Let $\{ f_1,\ldots,f_s \}$ be a set of binomials in $I_M$. Set $J=(f_1,\ldots,f_s)$ and $G=<{\bf y}_{\widehat{f_1}},\ldots,{\bf y}_{\widehat{f_s}}> \subset N(ker_{\bz}(M))$. If $char(K)=p \neq 0$ (resp. $char(K)=0$), then $I_M=rad(I_N+J)$ if and only if:\\
\noindent (a) $p^k N(ker_{\bz}(M)) \subset G$ for some $k \in \bn$ (resp. $N(ker_{\bz}(M))=G$),\\
\noindent (b) $rad((\phi(I_M))^{e},t_i)=rad(\phi(f_1),\ldots,\phi(f_s),t_i)$ for all $i=1,\ldots,l$.
\end{theorem}
\demo Suppose that $I_M=rad(I_N+J)$, then from Theorem 3.1 $rad((\phi(I_M))^{e})=rad(\phi(f_1),\ldots,\phi(f_s))$. Clearly $$rad((\phi(I_M))^{e},t_i)=rad(\phi(f_1),\ldots,\phi(f_s),t_i)$$ for all $i=1,\ldots,l$. In addition $rad((\phi(I_M))^{e}):(t_1 \cdots t_l)^{\infty}=rad(\phi(f_1),\ldots,\phi(f_s)):(t_1 \cdots t_l)^{\infty}$ and therefore $$rad((\phi(I_M))^{e}:(t_1 \cdots t_l)^{\infty})=rad((\phi(f_1),\ldots,\phi(f_s)):(t_1 \cdots t_l)^{\infty}).$$ So, from Lemma 4.1, we obtain the equality $rad(I_{N(ker_{\bz}(M))})=rad(I_G)$. Now Proposition 4.2 assures that in characteristic zero $N(ker_{\bz}(M))=G$, and in positive characteristic $p^k N(ker_{\bz}(M)) \subset G$ for some $k \in \bn$. Conversely suppose that (a) and (b) hold. By Proposition 4.2 $rad(I_{N(ker_{\bz}(M))})=rad(I_G)$, which implies that $rad((\phi(I_M))^{e}:(t_1 \cdots t_l)^{\infty})=rad((\phi(f_1),\ldots,\phi(f_s)):(t_1 \cdots t_l)^{\infty})$. If $I$ is any ideal of $K[t_1,\ldots,t_l]$, then by Lemma 3.2 in \cite{E-S} the radical of $I$ satisfies $$rad(I)=rad(I:(t_1 \cdots t_l)^{\infty}) \cap rad(I,t_1) \cap \cdots \cap rad(I,t_l).$$ Applying this formula to $(\phi(I_M))^{e}$ we obtain that $$rad((\phi(I_M))^{e})=rad(\phi(f_1),\ldots,\phi(f_s)).$$ Hence $I_M=rad(I_{N}+J)$. $\diamondsuit $
\begin{remark} {\rm In the special case that $I_N=(0)$ we take Theorem 2.5 in \cite{E-V}.}
\end{remark}
\begin{example} {\rm In \cite{Eli} S. Eliahou studied the binomial generation of the radical of the ideal of a monomial curve. Our theory will provide different binomial generators arising from different projections. For example, let $a \geq 7$ be an odd integer and let $M_a=(4,6,a,a+2)$. The toric ideal $I_{M_a}$ is a projection of the toric ideal $I_{D_a}$ associated with the matrix $$D_a=\left(\matrix{a-2 & a-4 & 2 & 0 \cr
0 & 2 & a-4 & a-2 \cr} \right).$$ Note that $I_{D_a}=rad(x_2^{a-2}-x_1^{a-4}x_4^2,x_3^{a-2}-x_1^2x_4^{a-4},x_1x_4-x_2x_3)$. Set $f_1=x_1^{a+2}-x_4^4,f_2=x_4^2-x_1x_3^2$ and $G=<(a+2,-4)>$. We have $V((\phi(I_{M_a}))^{e},t_i) \cap {\bar{K}}^{2}=V(\phi(f_1),\phi(f_2),t_i) \cap {\bar{K}}^{2}=\{ {\bf 0} \}$, since $\phi(f_1)=t_1^{(a-2)(a+2)}-t_2^{4(a-2)}$ and therefore $t_1=0$ if and only if $t_2=0$. Thus $rad((\phi(I_{M_a}))^{e},t_i)=rad(\phi(f_1),\phi(f_2),t_i)$ for $i=1,2$. Let ${\bf u}=(u_1,\ldots,u_4) \in ker_{\bz}(M_a)$. Then ${\bf y}_{\bf u}=(u_1+u_2+u_3+u_4)(a+2,-4)$ and therefore $D_a(ker_{\bz}(M_a))=G$. Now Theorem 4.4 assures that $I_{M_a}=rad(x_2^{a-2}-x_1^{a-4}x_4^2,x_3^{a-2}-x_1^2x_4^{a-4},x_1x_4-x_2x_3,x_1^{a+2}-x_4^4,x_4^2-x_1x_3^2)$.\\
But also $I_{M_a}$ is a projection of the toric ideal $I_{N_a}=(x_1^{\frac{a+1}{2}}-x_3x_4,x_1^3-x_2^2)$ associated with the matrix $$N_a=\left(\matrix{2 & 3 & a+1 & 0 \cr
2 & 3 & 0 & a+1 \cr} \right).$$ Let $\psi:K[x_1,\ldots,x_4] \rightarrow K[t_1,t_2]$ be the $K$-algebra homomorphism with $I_{N_a}=ker(\psi)$. Set $g_1=x_3^{a+2}-x_4^{a}$, $g_2=x_1x_4-x_2x_3$ and $H=<(a+2,-a)>$. Using the same arguments as before we take $rad((\psi(I_{M_a}))^{e},t_i)=rad(\psi(g_1),\psi(g_2),t_i)$ for $i=1,2$. If ${\bf u}=(u_1,\ldots,u_4) \in ker_{\bz}(M_a)$, then ${\bf y}_{\bf u}=(2u_1+3u_2+\frac{a+1}{2}u_3+\frac{a+1}{2}u_4)(a+2,-a)$ and therefore $N_{a}(ker_{\bz}(M_a))=H$. Consequently $I_{M_a}=rad(x_1^{\frac{a+1}{2}}-x_3x_4,x_1^3-x_2^2,x_3^{a+2}-x_4^{a},x_1x_4-x_2x_3)$.}
\end{example} 
\section{Set-theoretic complete intersection}   
Let $I_N \subset K[x_1,\ldots,x_n]$ be a toric ideal of height $r_1 \geq 1$ associated with an $l \times n$ integer matrix $N$ with non-negative entries. Let $B$ be the set of columns of $N$. Suppose that $\Gamma(N)=V(I_N)$ in ${\bar{K}}^n$.\\
\par We consider a lattice $L=ker_{\bz}(D)$ in $\bz^l$, where $D$ is an $m \times l$ rational matrix such that the matrix $M=DN$ has integer entries. The last statement means that $I_M$ is a projection of $I_N$. Let $r_2$ be the height of $I_M$ and $\phi:K[x_1,\ldots,x_n] \rightarrow K[t_1,\ldots,t_{l}]$ the $K$-algebra homomorphism with $I_N=ker(\phi)$. 
\begin{lemma} The dimension of the ${\bq}$-vector space $ker(D) \cap {\bq}B$ equals the difference $r_2-r_1$.
\end{lemma}
\demo Let $\{{\bf u}_1,\ldots,{\bf u}_{r_1} \}$ be a basis of $ker(N)$ and $\{{\bf u}_1,\ldots,{\bf u}_{r_1},{\bf u}_{r_{1}+1},\ldots,{\bf u}_{r_2} \}$ a basis of $ker(M)$. For the sake of simplicity the symbol ${\bf y}_{i}$ will represent ${\bf y}_{{\bf u}_i}$. We will show that $$ker(D) \cap {\bq}B=\bq \{{\bf y}_{{r_1}+1},\ldots,{\bf y}_{r_2} \}.$$ Obviously $\bq \{{\bf y}_{{r_1}+1},\ldots,{\bf y}_{r_2} \} \subseteq ker(D) \cap {\bq}B$. Let ${\bf v} \in ker(D) \cap {\bq}B$, then ${\bf v}={\bf y}_{\bf z}$ for some vector ${\bf z} \in \bq^n$. The vector ${\bf z}$ belongs to $ker(M)$, since $M=DN$. Thus ${\bf z}=\sum_{i=1}^{r_2} \kappa_{i}{\bf u}_i$ for some rationals $\kappa_{1},\ldots,\kappa_{r_2}$. Consequently ${\bf v}=\sum_{i=r_{1}+1}^{r_2} \kappa_{i}{\bf y}_{i} \in \bq \{{\bf y}_{r_{1}+1},\ldots,{\bf y}_{r_2} \}$. It remains to show that the set $\{{\bf y}_{r_{1}+1},\ldots,{\bf y}_{r_2} \}$ is linearly independent. Every relation of the form $\sum_{i={r_1}+1}^{r_2} \kappa_{i}{\bf y}_i={\bf 0}$ implies that the vector $\sum_{i={r_1}+1}^{r_2} \kappa_i{\bf u}_i$ belongs to $ker(N)$, so there exist some $\lambda_i$ such that $\sum_{i={r_1}+1}^{r_2} \kappa_{i}{\bf u}_i=\sum_{i=1}^{r_1} \lambda_{i}{\bf u}_i$. But the set $\{ {\bf u}_1,\ldots,{\bf u}_{r_2} \}$ is linearly independent, so all the $k_i$ are equal to zero. $\diamondsuit $
\begin{remark} {\rm The rank of the lattice $L \cap {\bz}B$ is equal to the dimension of $ker(D) \cap {\bq}B$. Also $L \cap {\bz}B$ coincides with the lattice $N(ker_{\bz}(M))$, so $$ht(I_{N(ker_{\bz}(M))})=r_2-r_1.$$ }
\end{remark}
\begin{theorem} Set $s=r_2-r_1$. If there are polynomials $f_1,\ldots,f_s$ in $I_M$ such that $rad(I_{N(ker_{\bz}(M))})=rad(\phi(f_1),\ldots,\phi(f_s))$, then $I_M$ is set-theoretic complete intersection on $I_N$. 
\end{theorem}
\demo We have that $rad((\phi(I_M))^{e}) \subseteq rad(I_{N(ker_{\bz}(M))})$. From the assumption $$rad(I_{N(ker_{\bz}(M))})=rad(\phi(f_1),\ldots,\phi(f_s)),$$ so $rad((\phi(I_M))^{e})=rad(\phi(f_1),\ldots,\phi(f_s))$ and therefore $I_M=rad(I_N+(f_1,\ldots,f_s))$.$\diamondsuit $
\smallskip
\newline
\par Combining Corollary 3.5 and Theorem 5.3 we get the following Corollary.
\begin{Corollary} Set $s=r_2-r_1$. If $I_N$ is set-theoretic complete intersection and $rad(I_{N(ker_{\bz}(M))})=rad(\phi(f_1),\ldots,\phi(f_s))$ for some polynomials $f_1,\ldots,f_s$ in $I_M$, then $I_M$ is set-theoretic complete intersection.
\end{Corollary}
\begin{example} {\rm In this example we will use the previous results to prove that the toric ideal $I_M$ of height $4$ associated with the matrix $$M= \left(\matrix{7 & 0 & 0 & 5 & 4 & 5 & 2 \cr
0 & 7 & 0 & 3 & 1 & 0 & 0 \cr
0 & 0 & 7 & 0 & 0 & 4 & 3 \cr} \right)$$ is set-theoretic complete intersection. Let $I_N=(x_5^3-x_1x_4,x_7^5-x_3x_6^2) \subset K[x_1,\ldots,x_7]$ be the toric ideal of height $2$ associated to the matrix $$N=\left(\matrix{3 & 0 & 0 & 0 & 1 & 0 & 0 \cr
0 & 3 & 0 & 0 & 0 & 0 & 0 \cr
0 & 0 & 5 & 0 & 0 & 0 & 1 \cr
0 & 0 & 0 & 3 & 1 & 0 & 0 \cr
0 & 0 & 0 & 0 & 0 & 5 & 2 \cr} \right).$$ Note that $I_M$ is a projection of $I_N$. Suppose that $\phi:K[x_1,\ldots,x_7] \rightarrow K[t_1,\ldots,t_5]$ is the $K$-algebra homomorphism with $I_N=ker(\phi)$. The set of vectors $$\{(1,0,0,1,-3,0,0),(0,0,1,0,0,2,-5) \}$$ constitutes a base for $ker_{\bz}(N)$ and the set $$\{(1,0,0,1,-3,0,0),(0,0,1,0,0,2,-5),(1,0,0,0,0,-3,4),(0,-1,0,4,-5,0,0) \}$$ constitutes a base for $ker_{\bz}(M)$. Thus $$N(ker_{\bz}(M))=<(3,0,4,0,-7),(-5,-3,0,7,0)>$$ and therefore $I_{N(ker_{\bz}(M))}=(t_5^7-t_1^3t_3^4,t_4^7-t_1^5t_2^3)$. Let $f_1=x_4^7-3x_1x_2x_4^4x_5^2+3x_1^3x_2^2x_4^2x_5-x_1^5x_2^3 \in I_M$ and $f_2=x_6^7-5x_1x_6^4x_7^4+10x_1^2x_3x_6^3x_7^3-10x_1^3x_3^2x_6^2x_7^2+5x_1^4x_3^3x_6x_7-x_1^5x_3^4 \in I_M$. We have $(t_4^7-t_1^5t_2^3)^3=\phi(f_1)$ and $(t_5^7-t_1^3t_3^4)^5=\phi(f_2)$. So $rad(I_{N(ker_{\bz}(M))})=rad(\phi(f_1),\phi(f_2))$, which implies that $I_M$ is the set-theoretic complete intersection of $x_5^3-x_1x_4,x_7^5-x_3x_6^2,x_4^7-3x_1x_2x_4^4x_5^2+3x_1^3x_2^2x_4^2x_5-x_1^5x_2^3,x_6^7-5x_1x_6^4x_7^4+10x_1^2x_3x_6^3x_7^3-10x_1^3x_3^2x_6^2x_7^2+5x_1^4x_3^3x_6x_7-x_1^5x_3^4$. } \\
\end{example}
\par We will compute the ideal $I_{N(ker_{\bz}(M))}$ in the special case that $s=1$. Let $$N=\pmatrix{b_1 & 0 & \ldots & 0 & d_{1,1} & \ldots & d_{r_{1},1} \cr
0 & b_2 & \ldots & 0 & d_{1,2} & \ldots & d_{r_{1},2} \cr
\ldots & \ldots & \ldots & \ldots & \ldots & \ldots & \ldots \cr
0 & 0 & \ldots & b_l & d_{1,l} & \ldots & d_{r_{1},l} \cr},$$
where $b_1,\ldots,b_l$ are positive integers and $d_{i,j}$ are non-negative integers such that, for all $i=1,\ldots,r_1$, at least one of $d_{i,1},\ldots,d_{i,l}$ is non zero. From Corollary 2 in \cite{K-T} we have $\Gamma(N)=V(I_N)$ in ${\bar{K}}^n$. The symbol $\left| N \right|$ will represent the greatest common divisor of the subdeterminants of $N$ of order $l$. Set $L=<{\bf a}>$ and $w=\frac{\left| N \right|}{\left| (N{\bf a}^{\bf T}) \right|}$, where $(N{\bf a}^{\bf T})$ is the augmented matrix. Given a vector ${\bf u}$ in $\bz^l$, the binomial ${\bf t}^{{\bf u}_+}-{\bf t}^{{\bf u}_-}$ will be denoted by $F({\bf u})$.
\begin{theorem} The lattice ideal $I_{N(ker_{\bz}(M))}$ is equal to the ideal generated by $F(w{\bf a})$. Moreover, if $I_N$ is set-theoretic complete intersection and there exists $g \in I_M$ such that $rad(\phi(g))=rad(F(w{\bf a}))$, then $I_M$ is set-theoretic complete intersection.
\end{theorem}
\demo We have $L \cap {\bz}B=<w{\bf a}>$, since every system of the form ${\bf y}_{\bf u}=c{\bf a}$ has a solution if and only if $c$ is an integer multiple of $w$ (\cite{BMT}, Theorem 1).Thus $N(ker_{\bz}(M))=<w{\bf a}>$ and therefore $I_{N(ker_{\bz}(M))}=(F(w{\bf a}))$. Also, from Corollary 5.4, $I_M$ is set-theoretic complete intersection. $\diamondsuit $
\begin{remark} {\rm In the example 4.6 the choice of $f_2=x_4^2-x_1x_3^2$ (resp. $g_2=x_1x_4-x_2x_3$) was made by solving the system ${\bf y}_{\bf u}=(a+2,-4)$ (resp. ${\bf y}_{\bf u}=(a+2,-a)$). }
\end{remark}  
\section{Applications} In this section we will present some applications of the theory developed in section 5.\\
\par We consider the toric ideal of height $d-1$ associated with the $(m+1) \times (m+d)$ integer matrix $$N_{d}=\left(\matrix{d & d-1 & d-2 & \ldots & 1 & 0 & 0 & \ldots & 0 \cr
0 & 1 & 2 & \ldots & d-1 & d & 0 & \ldots & 0 \cr
0 & 0 & 0 & \ldots & 0 & 0 & d & \ldots & 0 & \cr
\ldots & \ldots & \ldots & \ldots & \ldots & \ldots &\ldots & \ldots & \ldots \cr
0 & 0 & 0 & \ldots & 0 & 0 & 0 & \ldots & d \cr} \right),$$ where $d>1$. The toric ideal $I_{N_d} \subset K[x_1,\ldots,x_{m+d}]$ is set-theoretic complete intersection, for details see \cite{R-V}, \cite{Tho}. Theorem 6.1 will generalize this result. Let $c_1$ be a positive integer, $c_2,\ldots,c_{m+1}$ be non-negative integers with $c_1 \geq dc_2$ and let $L=ker_{\bz}(D)$ for $$D=\left(\matrix{\frac{c_1}{d} & \frac{c_1-dc_2}{d} & 0 & \ldots & 0 \cr
0 & c_3 & \frac{c_1}{d} & \ldots & 0 \cr
\ldots & \ldots & \ldots & \ldots & \ldots \cr
0 & c_{m+1} & 0 & \ldots & \frac{c_1}{d} \cr} \right).$$ Set ${\bf g}=\gcd(c_1,c_1-dc_2,dc_3,\ldots,dc_{m+1})$ and let ${c_1}^{\star}={\frac{c_1}{\bf g}}$, $(c_1-dc_2)^{\star}={\frac{c_1-dc_2}{\bf g}}$ and $(dc_i)^{\star}={\frac{dc_i}{\bf g}}$ for $i=3,\ldots,{m+1}$. Let ${\bf a}=((c_1-dc_2)^{\star},-{c_1}^{\star},(dc_3)^{\star},\ldots,(dc_{m+1})^{\star})$. Observe that $L=<{\bf a}>$.
\begin{theorem} The toric ideal of height $d$ associated with the $ m \times (m+d)$ matrix $$M_{c_1,\ldots,c_{m+1},d}=\left(\matrix{c_1 & c_1-c_2 & c_1-2c_2 & \ldots & c_1-dc_2 & 0 & \ldots & 0 \cr
0 & c_3 & 2c_3 & \ldots & dc_3 & c_1 & \ldots & 0 \cr
\ldots & \ldots & \ldots \ldots & \ldots & \ldots & \ldots & \ldots& \ldots \cr
0 & c_{m+1} & 2c_{m+1} & \ldots & dc_{m+1} & 0 & \ldots & c_1 \cr} \right)$$ is set-theoretic complete intersection.
\end{theorem}
\demo Without loss of generality we can assume that the greatest common divisor of the elements of $M_{c_1,\ldots,c_{m+1},d}$ is equal to $1$. Let $\phi:K[x_1,\ldots,x_{m+d}] \rightarrow K[t_1,\ldots,t_{m+1}]$ be the $K$-algebra homomorphism with $I_{N_d}=ker(\phi)$. Note that $I_{M_{c_1,\ldots,c_{m+1},d}}$ is a projection of $I_{N_d}$. For the integer ${\bf g}$ we have ${\bf g}/d$, since ${\bf g}/dc_i$ for all $i=2,\ldots,{m+1}$ and $\gcd({\bf g},c_2,\ldots,c_{m+1})=1$. In this case $\left| N_d \right| =d^m$, $\left| (N_d{\bf a}^{\bf T}) \right| =\frac{d^{m}}{\bf g}$ and $w={\bf g}$. Also $F(w{\bf a})=t_2^{c_1}-t_1^{c_1-dc_2}t_3^{dc_3} \cdots t_{m+1}^{dc_{m+1}}$. We have $(t_2^{c_1}-t_1^{c_1-dc_2}t_3^{dc_3} \cdots t_{m+1}^{dc_{m+1}})^d=\phi(f)$ for $$f=x_{d+1}^{c_1}-\pmatrix{d \cr 1 \cr}x_d^{c_1-dc_2}x_{d+1}^{(d-1)c_2}x_{d+2}^{c_3} \cdots x_{d+m}^{c_{m+1}}+\pmatrix{d \cr 2 \cr}x_{d-1}^{c_1-dc_2}x_{d+1}^{(d-2)c_2}x_{d+2}^{2c_3} \cdots x_{d+m}^{2c_{m+1}}-$$ $$\cdots + (-1)^{d-1}\pmatrix{d \cr {d-1} \cr}x_2^{c_1-dc_2}x_{d+1}^{c_2}x_{d+2}^{(d-1)c_3} \cdots x_{d+m}^{(d-1)c_{m+1}} + (-1)^{d}x_1^{c_1-dc_2}x_{d+2}^{dc_3} \cdots x_{d+m}^{dc_{m+1}}.$$ Notice that $f$ belongs to $I_{M_{c_1,\ldots,c_{m+1},d}}$. Consequently, from Theorem 5.6, the toric ideal $I_{M_{c_1,\ldots,c_{m+1},d}}$ is set-theoretic complete intersection. $\diamondsuit $      
\smallskip
\newline
\par Next we prove that the toric ideal associated to the row matrix $M_{a,b}=(a,a+2b,2a+3b,2a+5b)$ is set-theoretic complete intersection. Especially, when $a=4, b=1$ we deduce that the ideal of the monomial curve $(t^4,t^6,t^{11},t^{13})$ is set-theoretic complete intersection. We consider the toric ideal associated to the matrix $$N=\left(\matrix{5 & 1 & 4 & 0 \cr
0 & 2 & 3 & 5 \cr} \right).$$ The toric ideal $I_N \subset K[x_1,\ldots,x_4]$ is the set-theoretic complete intersection of $x_3^2-x_1x_2^3$ and $x_2^5-2x_2x_3x_4+x_1x_4^2$. Let $L=ker_{\bz}(E)$ for $E=(\frac{a}{5},\frac{2a+5b}{5})$. Set ${\bf h}=\gcd(a,2a+5b)$, $a^{\star}={\frac{a}{\bf h}}$ and $(2a+5b)^{\star}={\frac{2a+5b}{\bf h}}$. Note that $L=<{\bf a}>$, where ${\bf a}=(-(2a+5b)^{\star},a^{\star})$. In addition $I_{M_{a,b}}$ is a projection of $I_{N}$.
\begin{theorem} For every positive integers $a,b$ the ideal of the monomial curve $$(t^{a},t^{a+2b},t^{2a+3b},t^{2a+5b})$$ is set-theoretic complete intersection.
\end{theorem}
\demo Let $\phi:K[x_1,x_2,x_3,x_4] \rightarrow K[t_1,t_2]$ be the $K$-algebra homomorphism with $I_N=ker(\phi)$. For the integer ${\bf h}$ we have ${\bf h}/a$ and ${\bf h}/5b$, so ${\bf h}/5$ since $gcd({\bf h},b)=1$. Here $\left| N \right| =5$, $\left| (N{\bf a}^{\bf T}) \right| =\frac{5}{\bf h}$ and $w={\bf h}$. Also $F(w{\bf a})=t_1^{2a+5b}-t_2^{a}$. When $a=1$ the toric ideal associated to the matrix $M_{1,b}$ is obviously set-theoretic complete intersection. Suppose that $a>1$, which implies that $a=2\mu+3\nu$ for some non-negative integers $\mu, \nu$. We have $(t_1^{2a+5b}-t_2^a)^5=\phi(f)$, where $f=x_1^{4\mu+6\nu+5b}-5x_1^{3\mu+4\nu+4b}x_2^{\mu}x_3^{\nu}+10x_1^{2\mu+2\nu+3b}x_2^{2\mu}x_3^{2\nu}-10x_1^{\mu+2b}x_2^{3\mu}x_3^{3\nu}+5x_1^{\nu+b}x_2^{\nu}x_3^{\mu}x_4^{\mu+2\nu}-x_4^{2\mu+3\nu}$ belongs to $I_{M_{a,b}}$. Therefore, from Theorem 5.6, the ideal of the curve $(t^{a},t^{a+2b},t^{2a+3b},t^{2a+5b})$ is set-theoretic complete intersection. $\diamondsuit $
\smallskip
\newline
\par Finally we prove that the toric ideal associated with the matrix $M_{a,b}=(4,6,a,b)$ is set-theoretic complete intersection.
\begin{theorem} For every positive integers $a,b$ the ideal of the monomial curve $(t^4,t^6,t^{a},t^{b})$ is set-theoretic complete intersection.
\end{theorem}
\demo Suppose that $b=a+k$, where $k$ is a positive integer. If $a$ or $b$ is even, then the semigroup $\bn(4,6,a,b)$ is symmetric by Proposition 2.1 in \cite{He} and therefore $I_{M_{a,b}}$ is set-theoretic complete intersection (see \cite{Bre}). It remains to examine the case $a$ is odd and $k$ is even. When $k \geq 4$, the semigroup $\bn (4,6,a,b)$ is symmetric and the result is straightforward. Therefore we have to deal only with the case $k=2$. Since $a>1$, there is a non-negative integer $\mu$ and a positive integer $\nu$ such that $a=2 \mu  +3 \nu$. We consider the toric ideal $I_{N_a}$ associated with the matrix $N_a$ of the example 4.6. In this example it was proved that $N_{a}(ker_{\bz}(M_{a,a+2}))=<(a+2,-a)>$ and so $I_{N_{a}(ker_{\bz}(M_{a,a+2}))}=(t_1^{a+2}-t_2^{a})$. Set $$f=\sum_{0 \leq i \leq \frac{a+1}{2}} (-1)^{i}\pmatrix{a+1 \cr i \cr}x_1^{i \mu}x_2^{i \nu}x_3^{a+2-2i}+\sum_{\frac{a+1}{2}<i \leq {a+1}} (-1)^{i}\pmatrix{a+1 \cr i \cr}x_1^{(a+1-i)(\mu+1)}x_2^{\nu(a+1-i)}x_4^{2i-a-2}$$ and observe that $f \in I_{M_{a,a+2}}$. We have $(t_1^{a+2}-t_2^{a})^{a+1}=\psi(f)$ and therefore $I_{M_{a,a+2}}$ is the set-theoretic complete intersection of $f,x_1^{\frac{a+1}{2}}-x_3x_4,x_1^3-x_2^2$. $\diamondsuit $       
\begin{remark} {\rm The last Theorem provides the polynomials that minimally generate up to radical the ideal of the Eliahou's curve $(t^4,t^6,t^7,t^9)$, see also \cite{El} for a proof that the above ideal is set-theoretic complete intersection. Also it provides a different minimal polynomial generating set, than the one obtained in Theorem 6.2, up to radical for the ideal of the monomial curve $(t^4,t^6,t^{11},t^{13})$.}
\end{remark}

\par Department of Mathematics,
Section of Algebra and Geometry, University of Ioannina, Ioannina 45110 (GREECE) \\
\par Email address: akatsabekis@in.gr \\
\end{document}